\documentclass[10pt]{amsart}
\usepackage{graphicx}
\usepackage{psfrag}
\usepackage{amssymb}
\def\op{\overline\partial}

\let\newpf\proof \let\proof\relax 
\newenvironment{pf}{\newpf[\proofname]}{\qed\endtrivlist}

\def\cal{\mathcal}

\newtheorem{thm}{Theorem}[section]
\newtheorem{problem}{Problem}
\newtheorem{cor}[thm]{Corollary}

\newtheorem{lemma}[thm]{Lemma}

\theoremstyle{remark}
\newtheorem{rem}{Remark}[section]

\numberwithin{equation}{section}

\def \bn {\hfill \\ \smallskip\noindent}

\theoremstyle{definition}

\def\proof{\bn {\bf Proof.} }

\def\note#1
{\marginpar
{\tiny $\leftarrow$
\par
\hfuzz=20pt \hbadness=9000 \hyphenpenalty=-100 \exhyphenpenalty=-100
\pretolerance=-1 \tolerance=9999 \doublehyphendemerits=-100000
\finalhyphendemerits=-100000 \baselineskip=6pt
#1}\hfuzz=1pt}

\newcommand{\Rat}{\operatorname{Rat}}

\newcommand{\orb}{\operatorname{orb}}

\newcommand{\id}{\operatorname{id}}

\newcommand{\Dil}{\operatorname{Dil}}

\newcommand{\esssup}{\operatorname{ess-sup}}

\renewcommand{\epsilon}{{\varepsilon}}

\newcommand{\CC}{{\cal C}}

\newcommand{\HH}{{\cal H}}

\newcommand{\MM}{{\cal M}}

\newcommand{\C}{{\mathbb C}}
\newcommand{\D}{{\mathbb D}}

\def\B0{{\bold{0}}}

\catcode`\@=12

\def\Empty{}
\newcommand\oplabel[1]{
  \def\OpArg{#1} \ifx \OpArg\Empty {} \else
  	\label{#1}
  \fi}
		
\newcommand{\comm}[1]{}
\newcommand{\comment}[1]{}

\begin{document}

\bigskip\bigskip

\title{Infinitesimal perturbations of rational maps}
\author {Artur Avila}


\address{
Coll\`ege de France, 3 Rue d'Ulm, 75005 Paris -- France.
}

\email{avila@impa.br}

\date{\today}

\begin{abstract} 

We analyze the infinitesimal effect of holomorphic perturbations of the
dynamics of a structurally stable
rational map on a neighborhood of its Julia set.  This implies some
restrictions on the behavior of critical points.

\end{abstract}

\setcounter{tocdepth}{1}

\maketitle

\section{Introduction}

A major open problem in one-dimensional dynamics is to show that
hyperbolic rational maps are dense among all rational maps (Fatou
conjecture).  After the work of Ma\~n\'e-Sad-Sullivan \cite {MSS},
this problem reduces to show that all
structurally stable rational maps are hyperbolic.  In other words,
one should show that non-hyperbolic rational
maps can be perturbed to a map with different dynamics (in the topological
sense).

An obvious difficulty to perturb
rational maps is that there are not enough perturbations, since the space of
rational maps of a given degree is finite dimensional.
It is then hard to show that any of
them has non-trivial effects on the dynamics.

A natural approach is then to perturb rational maps inside a broader
class of dynamical systems, and later realize the perturbation as a rational
one.  One such class consists of quasiregular maps (the quasiconformal
analogous to rational maps).  The space of
quasiregular perturbations of a given rational map is infinite dimensional,
and it is indeed much simpler to change the dynamics of a rational map
through a quasiregular perturbation.  Of course one is then left with the
problem of realizing the perturbation as a rational one: this problem was
attacked by Thurston in some simple situations, and is the basic
problem in \cite {R}, where it is solved under certain geometric
assumptions on the recurrence of critical points.
Another application of quasiregular perturbations,
due to Shishikura \cite {Sh}, is a sharp bound on the number of
non-repelling periodic orbits of rational maps of a given degree.

(One should note the ressemblence to a usual approach to rigidity of
rational maps: to show that two maps are affinely conjugate, one first looks
for qc conjugacies.)

Another approach to analyzing perturbations of rational maps is to observe
that the Julia set of a structurally stable rational map $R$ depends
holomorphically on a neighborhood of $R$.  Linearizing this holomorphic
dependence for the case of a infinitesimal perturbation $R+\lambda v$
we get a representation $v=\alpha \circ R-DR \alpha$ for some
$\alpha$ which is a qc vector field on $J(R)$.  Looking for non-trivial
perturbations is then reduced to finding some $v$ which can not satisfy such
an equation.  A successful application of those ideas is for instance
Ma\~n\'e's proof of instability of Herman rings \cite {Ma}. For other
applications, the fact that the space of such $v$ is finite dimensional
can be a major obstacle, and this is the problem we try to address here.

Our aim in this paper is to present a technique of infinitesimal analysis of
certain quasiregular perturbations.  We consider a rational map $R$ and a
quasiregular perturbation $R+\lambda v$ of $R$ which is holomorphic on a
neighborhood of the Julia set $J(R)$.
This is enough to increase a lot our
freedom if we assume $J(R) \neq \overline \C$.  We then show that if
$R$ is structurally stable then the ``Julia set'' of $R+\lambda v$
depends holomorphically on
$\lambda$ for small $\lambda$.  Linearizing this holomorphic dependence we
get the representation $v=\alpha \circ R-DR \alpha$ on $J(R)$.  Thus we
increased the space of ``test functions'' against which to check the
possibility of the above representation from a finite dimensional space to
a much bigger space.

We present a direct application of this
procedure by considering the holomorphic dependence of
iterates of critical points $c_\lambda$
of $R+\lambda v$.  The fact that the speed of motion of those iterates must
be bounded reflects on the fact that certain sequences
depending on $v$ are bounded, and this imposes restrictions
on the dynamics of $R$.

As a corollary of this analysis, we are able to show that a rational
map such that $J(R) \subset \C$ (the case $J(R) \neq \overline \C$
reduces to this one by a Moebius change of coordinates),
with one critical point satisfying a summability condition
and an additional technical condition
(which can be interpreted as `the space of
holomorphic vector fields in a neighborhood of $J(R)$ is big enough'')
is not structurally stable.  More precisely we have:

\begin{thm}[see Corollary \ref {nontrivial}] \label {thmnt}

Let $R$ be a rational map such that $J(R) \subset \C$ and such that
there exists a critical point $c \in J(R)$ with
\begin{equation} \label {sumc}
\sum_{k=0}^\infty \frac {1} {|DR^k(R(c))|}<\infty.
\end{equation}
If there exists a vector field $v$ which is
holomorphic on a neighborhood of the closure of the orbit of $c$
and such that
\begin{equation} \label {vhol}
\sum_{k=0}^\infty \frac {v(R^k(c))} {DR^k(R(c))} \neq 0
\end{equation}
then $R$ is not structurally stable.

\end{thm}

(Notice that the summability condition (\ref {sumc})
implies (\ref {vhol}) for many {\it continuous} vector fields $v$.)

This recovers a result of Levin
and Makienko (\cite {Le}, \cite {M}),
which is based on a different approach (analysis of a Ruelle transfer
operator).  Indeed this work was
motivated as an attempt to give a geometric interpretation to those results,
in the hope that both points of view could be eventually merged.

We finally discuss the meaning of this technique in the better understood
case of polynomial maps.

Although we do not develop it here, infinitesimal
analysis is also an important technique in real one-dimensional dynamics.
Indeed this path was thoroughly developed in the case of unimodal maps in
\cite {ALM}, and the results of that work could shed more light to
the setting we consider here.

{\bf Acknowledgements:}  The author thanks Michael Benedicks,
Juan Rivera-Letelier and Peter Makienko for many fruitful discussions
and to Lukas Geyer for pointing out a mistake in a preliminary version.

\subsection{Notation and definitions}

We let $\C$ denote the complex plane
and $\overline \C$ the Riemann Sphere.
We let $\D_\epsilon \subset \C$
denote the disk with radius $\epsilon$ around $0$.  The identity map is
denoted $\id$.

We assume the reader is familiar with the theory of quasiconformal maps. 
The Beltrami differential of a qc map $h$ is $\mu=\op h/\partial h$.  The
dilatation of $h$ is $\Dil(h)=\esssup (1+|\mu|)/(1-|\mu|)$.  $h$ is said to
be normalized if it fixes $0$, $1$ and $\infty$.  We will also say that
$h:X \to \overline \C$ is a qc map on $X$ if it extends to a (global)
qc map.

A quasiconformal vector field $\alpha$ of $\overline \C$ is a continuous
vector field with locally integrable distributional derivatives
$\op \alpha$ and $\partial \alpha$
in $L^1$ and $\op \alpha \in L^\infty$.  If $\alpha$ is only defined
on a set $X \subset \overline \C$ then we say $\alpha$ is qc on $X$
if it has a qc extension.

Let $U$ be a complex manifold with a base point $\lambda_0 \in U$.  Let $X
\subset \overline \C$ be any set.
A holomorphic motion of $X$ over $U$ (based on  $\lambda_0$) is a family
$H_\lambda:X \to \overline \C$, $\lambda \in U$, of injective
maps $H_\lambda$ such that for each $x \in X$,
$\lambda \mapsto H_\lambda(x)$ is holomorphic and $H_{\lambda_0}=\id$.
In this case each
$H_\lambda$ is qc on $X$ and the continuous extension
$H_\lambda:\overline X \to \C$ is also a
holomorphic motion.

In this paper we will usually let $U \subset \C$ be a domain
containing $0$ and consider holomorphic motions based on $0$.  We observe 
that the Beltrami differential $\op H_\lambda/\partial H_\lambda$ of a
holomorphic motion depends holomorphically on $\lambda$, and
$$
\alpha=\left .\frac {dH_\lambda} {d\lambda} \right |_{\lambda =0}
$$
is a qc vector field on $X$ (see \cite {ALM}).

\section{Rational maps}

Let us quickly recall the basic theory of rational maps (see \cite {MS}).
Let $\Rat_d$ be the space of rational maps of degree $d$.
Let $R \in \Rat_d$.
We define as usual the Fatou set of $R$, $F(R)$, as the biggest open set
where $\{R^n\}$ is a normal family.  The Julia set of $R$ is then defined as
$J(R)=\overline \C \setminus F(R)$.  The Julia set is the closure of
repelling periodic orbits of $R$, and there are at most finitely many
periodic orbits in $F(R)$.  It is easy to see that
$R$ takes connected components of $F(R)$ to connected
components of $F(R)$.  Moreover, according to Sullivan, every connected
component $U \subset F(R)$ is either periodic ($R^k(U)=U$ for some $k>0$) or
is eventually mapped to some periodic component, and there is only a finite
number of such periodic components.

We say that $R_0$ is structurally stable if there exists a neighborhood $U
\subset \Rat_d$ of $R_0$ such that for every $R \in U$,
there exists a homeomorphism $h_R:\overline \C \to \overline \C$ such that
$h_R \circ R_0 \circ h_R^{-1}=R$.  In this case the $h_R$ can be chosen to
form a holomorphic motion over $U$ based on $R_0$ (see \cite {MS}, \S 7).
In particular, $\lim_{R \to R_0}
\Dil h_R=1$ and $\lim_{R \to R_0} h_R=\id$ uniformly on $\C$.

It is known that if $R$ is structurally stable then every periodic
component $U$ of $F(R)$ (of period $k$)
must be attracting, that is, there exists $w \in U$ such that $R^k(w)=w$,
$|DR^k(w)|<1$ and for every $z
\in U$, $f^{kn}(z) \to w$.

Let $\Rat^\infty_d \subset \Rat_d$ be the space of all
rational maps with degree $d$ such that $\infty \in F(R)$.
Let $\Rat^a_d \subset \Rat_d$ be the space of maps for which
all periodic components of $F(R)$ are attracting and $Rat^s_d$ be the space
of structurally stable maps $R$.

\subsection{Quasiregular maps}

A map $\tilde R:\overline \C \to \overline \C$ is called a quasiregular map
if it is the composition of a rational map $R$ with a
qc map $h$, $\tilde R=R \circ h$.  One defines Beltrami differentials
and dilatation of quasiregular maps in the natural way.

Let us now consider a map $\tilde R$ which is a composition
of a qc map $h$ and a rational map $R$, $\tilde R=h \circ R$.  Let
$\tilde h$ be a qc map with the same Beltrami differential of $\tilde R$.
Then
$\tilde R \circ \tilde h^{-1}$ is a rational map.  So $\tilde R$
can also be seen as a composition of a rational map and a qc map
$\tilde R=(\tilde R \circ \tilde h^{-1}) \circ \tilde h$, and is a
quasiregular map.  It follows:

\begin{lemma}

The composition of quasiregular maps is quasiregular, furthermore
$\Dil(\tilde R_1 \circ \tilde R_2) \leq \Dil(\tilde R_1) \Dil(\tilde R_2)$.

\end{lemma}

Given a quasiregular map $\tilde R$, let us consider the support $K$ of
its Beltrami differential $\op \tilde R/\partial \tilde R$.  Let us say
that $\tilde R$ is tame if there exists $n>0$ such that for $m \geq n$,
$\tilde R^m(K) \cap K=\emptyset$.  Notice that in this case, for all
$m>0$,
$\Dil (\tilde R^m) \leq \Dil (\tilde R)^n$.

\begin{lemma} \label {bounded}

Let $\tilde R$ be a tame quasiregular map.
Then there exists a normalized qc map $h$ such that
$\Dil(h) \leq \sup \Dil(\tilde R^n)$ and $h \circ \tilde R \circ
h^{-1}$ is a rational map.

\end{lemma}

\begin{pf}

Let $\mu_n$ be the Beltrami differential of $\tilde R^n$, so that
$\mu_{n+1}=\tilde R^*(\mu_n)$.  Since $\tilde R$ is tame, it follows that 
for every $x$, the sequence $\mu_n(x)$ is eventually constant.  Let
$\mu$ be the pointwise limit of $\mu_n$ and let $h$ be a normalized
quasiconformal map with Beltrami differential $\mu$.  We have
$\tilde R^*(\mu)=\mu$, so that
$h \circ \tilde R \circ h^{-1}$ is a rational map.
\end{pf}

\begin{rem}

The conclusion of the above Lemma also holds if the tame assumption is
replaced by the weaker condition $\sup \Dil (\tilde R^n)<\infty$.

\end{rem}

\subsection{Quasiregular perturbations}

An admissible perturbation of a rational map $R$ is a family $R_\lambda$,
$\lambda \in \D_\epsilon$, of quasiregular maps such that

\begin{itemize}

\item $R_0=R$,

\item $\lim_{\lambda \to 0} R_\lambda(z)=R(z)$ uniformly on $\overline \C$,

\item $\lim_{\lambda \to 0} \Dil(R_\lambda)=1$,

\item there is a neighborhood $V_1$ of $J(R)$ such that $(\lambda,z) \mapsto
R_\lambda(z)$ is holomorphic on $\D_\epsilon \times V_1$,

\item there is a neighborhood $V_2$ of the attracting periodic orbits of $R$
and of the critical points of $R$ on $F(R)$ such that $R_\lambda=R$.

\end{itemize}

\begin{lemma} \label {qrpert}

Let $R_\lambda$ be an admissible family through $R \in \Rat^a_d$.  Then
there exists $n>0$, $\delta>0$ and a family of qc maps $h_\lambda$,
$\lambda \in \D_\delta$ such that
$h_\lambda^{-1} \circ R_\lambda \circ h_\lambda \in Rat_d$,
$\lim_{\lambda \to 0} h_\lambda=\id$ and
$\Dil(h_\lambda) \leq \Dil(R_\lambda)^n$.

\end{lemma}

\begin{pf}

Let $V_1$ and $V_2$ be as in the definition of admissible family and let
$K=\overline \C \setminus V_1 \setminus V_2$, so that for any $\lambda$,
$\op R_\lambda=0$ out of $K$.  Let $V \subset V_2$
be a neighborhood of the attracting cycles of
$R$ with $R(\overline V) \subset V$.
Since $R \in \Rat^a_d$, there exists $n>0$ such that
$R^n(K) \subset V$.  Taking $\delta$ small we may assume that
$R_\lambda^n(K) \subset V$ for $|\lambda|<\delta$.
By Lemma \ref {bounded}, there exists a family $h_\lambda$ of normalized 
qc maps such that $h_\lambda^{-1} \circ R_\lambda \circ h_\lambda \in
\Rat_d$, and $\Dil h_\lambda \leq (\Dil R_\lambda)^n$.  In particular,
$\lim_{\lambda \to 0} \Dil h_\lambda=1$ and since the $h_\lambda$ are
normalized, $\lim_{\lambda \to 0} h_\lambda=\id$ in the uniform topology. 
\end{pf}

\begin{rem} \label {qrpertrem}

The $h_\lambda$ we constructed actually satisfies $\op h_\lambda=0$ on
$J(R)$.  Indeed,
it is clear that $h_\lambda(J(R)) \subset V_1$ for $\lambda$ small.  In
particular, $R_\lambda^n$ is holomorphic on
$h_\lambda(J(R))$ for all $n$.  It is
easy to check that the proof of Lemma \ref {bounded} gives $\op h_\lambda=0$
on $J(R)$.

\end{rem}

\begin{lemma}

Let $R_\lambda$ be an admissible family through $R \in \Rat^s_d$.  Then
there exists a family of qc maps $H_\lambda$,
$\lambda \in \D_\delta$ such that
$H_\lambda^{-1} \circ R_\lambda \circ H_\lambda=R$,
$\lim_{\lambda \to 0} H_\lambda=\id$ and
$\lim_{\lambda \to 0} \Dil(H_\lambda)=1$.

\end{lemma}

\begin{pf}

The previous lemma gives us a family
$\tilde R_\lambda=h_\lambda^{-1} \circ R_\lambda
\circ h_\lambda \in Rat_d$ which is continuous on $\lambda=0$ (since
$\lim_{\lambda \to 0} h_\lambda=\id$).  Since $R$ is structurally stable,
there exists a family of quasiconformal maps $\tilde h_\lambda$ such that
$\tilde h_\lambda^{-1} \circ \tilde R_\lambda \circ \tilde h_\lambda=R$,
and
moreover $\lim \tilde h_\lambda=\id$ and $\lim \Dil(\tilde
h_\lambda)=1$.  Setting $H_\lambda=h_\lambda \circ \tilde h_\lambda$ we
get
the result.
\end{pf}

\begin{thm}

Let $R_\lambda$ be an admissible family through $R \in \Rat^s_d$, and let
$H_\lambda$ be a family of qc maps over $\D_\delta$ as above
such that $H_\lambda \circ R \circ H_\lambda^{-1}=R_\lambda$.
Then there exists $\epsilon>0$ and a holomorphic
motion $h_\lambda$ of $J(R)$ over $\D_\epsilon$
such that $h_\lambda(J(R))=H_\lambda(J(R))$ and $h_\lambda \circ R \circ
h_\lambda^{-1}=R_\lambda$.

\end{thm}

\begin{pf}

Let $W$ be a neighborhood of $J(R)$ where $R_\lambda$ is holomorphic and
such that all periodic orbits of $R$ in $W$ in fact belong to $J(R)$.
Let $J(R) \subset W' \subset \overline W' \subset W$ be a 
smaller neighborhood.
Taking $\epsilon<\delta$ small, we may assume that $J(R) \subset
H_\lambda^{-1}(W') \subset \overline {H_\lambda^{-1}(W')}
\subset W$ for all $\lambda \in \D_\epsilon$.

It is easy to see that $p$ is a periodic orbit of $R$ on $H_\lambda^{-1}(W')$
(of period $n$ and multiplier $\rho=DR^n(p)$) if and only
if $H_\lambda(p) \in W'$ is a periodic orbit of $R_\lambda$
(with period $n$ and multiplier $\rho_\lambda$) and
$\ln(|\rho|)/K \leq \ln(|\rho_\lambda|) \leq K \ln(\rho)$
with $K$ only depending only on $\Dil(H_\lambda)$, and in particular, since
the former is uniformly bounded, does not depend on $\lambda$.
So $H_\lambda(J(R))$ is the closure of the set of repelling
periodic orbits of $R_\lambda$ contained in $\overline W'$.

Let $P_n(\lambda)$ be the set
of repelling periodic orbits of period less than $n$ of $R_\lambda$ on
$H_\lambda(J(R))$.  The above discussion shows that
$H_\lambda(P_n(0))=P_n(\lambda)$.
Fixing $\lambda_0$, if $p_{\lambda_0} \in
P_n(\lambda_0)$ and $p_\lambda$, $\lambda$ near $\lambda_0$
is its holomorphic continuation (as a repelling periodic orbit), we have
$p_\lambda \in P_n(\lambda)$.  Since the $P_n(\lambda)$ are finite sets with
the same cardinality and $\D_\delta$ is simply connected,
it follows that there exists a unique holomorphic motion
$h_\lambda^n:P_n(0) \to P_n(\lambda)$, $\lambda \in \D_\delta$.
Passing to the limit $n \to \infty$ and
using the extension theorem for holomorphic motions of \cite {MSS}, we
obtain the desired holomorphic motion of $J(R)$.
\end{pf}

\begin{rem}

It is easy to see that if $c \in J(R)$ is a critical point of $R$ then
$h_\lambda(c)$ is a critical point of $R_\lambda$.

\end{rem}

\subsection{Infinitesimal perturbations}

Let now $R \in \Rat^\infty_d$.  Let us consider the space of
vector fields $v$ which are holomorphic on a neighborhood of $J(R)$.  This
space can be interpreted as a tangent space to certain quasiregular
perturbations due to the following obvious lemma.

\begin{lemma}

Let $R \in \Rat_d^\infty$ and let $v$ be a holomorphic vector field on a
neighborhood of $J(R)$.  Then there exists an admissible family
$R_\lambda$ through $R$ and a bounded neighborhood $V$
of $J(R)$ such that $R_\lambda=R+\lambda v$ on a neighborhood of $J(R)$.

\end{lemma}

\begin{pf}

Let $\tilde v$ be a $C^\infty$ vector field coinciding with $v$ on a
neighborhood of $J(R)$ and vanishing on a neighborhood of
the union of attracting periodic
orbits of $R$, critical points on $F(R)$, $R^{-1}(\infty)$ and $\infty$. 
Then it is easy to check that
$R_\lambda=R+\lambda \tilde v$ is quasiregular for
$\lambda$ small and has the required properties.
\end{pf}

\begin{cor}

Let $R \in \Rat^s_d \cap \Rat^\infty_d$ and let $v$ be a holomorphic vector
field on a bounded neighborhood $V$ of $J(R)$.
Then there exists $\delta>0$ and a
holomorphic motion $h_\lambda$, $\lambda \in \D_\delta$ of $J(R)$ such that
$h_\lambda(J(R)) \subset V$ and $h_\lambda \circ R \circ
h_\lambda^{-1}=R+\lambda v$.  Moreover if $c \in J(R)$ is a critical point
then $h_\lambda(c)$ is a critical point of $R+\lambda v$.

\end{cor}

\begin{thm}

Let $R \in \Rat^s_d \cap \Rat^\infty_d$ and let $v$ be a holomorphic vector
field on a neighborhood of $J(R)$.  Then there exists a unique
qc vector field $\alpha$ such that $v=\alpha \circ R-DR \alpha$ on $J(R)$.

\end{thm}

\begin{pf}

Just take
$$
\alpha=\left .\frac {d} {d\lambda} h_\lambda\right |_{\lambda=0}
$$
and linearize $h_\lambda \circ R \circ h_\lambda^{-1}=R+\lambda v$.
\end{pf}

\begin{cor} \label {struc}

Let $R \in \Rat^s_d \cap \Rat^\infty_d$ and let $v$ be a holomorphic vector
field on a neighborhood of $J(R)$.  Then for any $c \in J(R)$ critical point
of $R$, the sequence
$$
\left .\frac {d} {d\lambda} (R+\lambda v)^n(c) \right |_{\lambda=0}=
DR^{n-1}(R(c)) \sum_{k=0}^{n-1} \frac {v(R^k(c))} {DR^k(R(c))}
$$
is uniformly bounded.

\end{cor}

\begin{pf}

Let $\alpha$ be as in the previous Theorem.  Then it is easy to see by
induction that
$$
\alpha(R^n(c))=
\left .\frac {d} {d\lambda} (R+\lambda v)^n(c) \right |_{\lambda=0}=
DR^{n-1}(R(c)) \sum_{k=0}^{n-1} \frac {v(R^k(c))} {DR^k(R(c))}.
$$

(Since $R$ is structurally stable, there are no critical relations
so $DR^k(R(c))$ is non-zero for $k \geq 0$.)

Since $\alpha$ is a qc vector field, it must be uniformly bounded on $J(R)$,
so on $R^n(c)$.

(Another reasoning is that the above formulas give the speed of motion of
the $n$-th iterate of the critical point $c$ after perturbation by $v$. 
Since those iterates depend holomorphically over a definite neighborhood of
$0$, so Cauchy estimates give the desired bound.  This conclusion can be
obtained without the use of holomorphic motions.)
\end{pf}

\subsection{Measures}

Let us say that a critical point $c \in J(R)$ of $R$ is summable if
$$
\sum_{k=0}^\infty \frac {1} {|DR^k(R(c))|}<\infty.
$$
In this case, for a continuous map $v:\overline \C \to \C$ we let
$$
\mu_{R,c}(v)=\sum_{k=0}^\infty \frac {v(R^k(c))} {DR^k(R(c))},
$$
so that $\mu_{R,c}$ is a complex measure.
It is clear that $\mu_{R,c}$ does
not vanish over continuous functions.

Let us say that $\mu_{R,c}$ is non-trivial if it
does not vanish over the space of holomorphic functions on a
neighborhood of $\overline {\orb_R(c)}$.
It should be noted that if $\overline {\orb_R(c)}$ does not
disconnect the plane (for instance, if it is a Cantor set) then $\mu_{R,c}$
is non-trivial by the Mergelyan Approximation Theorem.

With those definitions, Theorem \ref {thmnt} can be restated as follows:

\begin{cor} \label {nontrivial}

Let $R \in \Rat_d^\infty$ and assume that
$c$ is a summable critical point such that
$\mu_{R,c}$ is non-trivial.
Then $R$ is not structurally stable.

\end{cor}

\begin{pf}

Assume by contradiction that $R$ is structurally stable, and let $v$ be a
holomorphic vector field on a neighborhood of $\overline {\orb_R(c)}$ such
that $\mu_{R,c}(v) \neq 0$.  By the Runge Approximation Theorem,
there exists a meromorphic vector field $\tilde v$ on
$\overline \C$ which is uniformly
close to $v$ on $\overline {\orb_R(c)}$.  Since $J(R) \neq \overline \C$,
it has empty interior so we may assume that the poles of
$\tilde v$ do not intersect $J(R)$, so $\tilde v$
is holomorphic in a neighborhood of $J(R)$.

Then
$$
\mu_{R,c}(\tilde v)=\lim \frac {1} {DR^{n-1}(R(c))}
\left .\frac {d} {d\lambda} (R+\lambda v)^n(c) \right |_{\lambda=0}.
$$
Since $c$ is a summable critical point, $|DR^n(R(c))| \to \infty$, so by
Corollary \ref {struc}, $\mu_{R,c}(\tilde v)=0$.
This is a contradiction since $\mu_{R,c}(\tilde v)$ approximates
$\mu_{R,c}(v)$.
\end{pf}

\begin{rem}

When $\omega(c)$ is a non-minimal set (this is often the case if
$\overline {\orb_R(c)}$ is not a Cantor set)
it is usually possible to tackle succesfully
the problem of structural stability using other techniques.

\end{rem}

\begin{rem}

Peter Makienko has informed me that
more general conditions for non-triviality can be obtained
using the theory of universal algebras.  For instance,
if $J(R)$ is the boundary of one of
the Fatou components (this is the case for polynomials), $\mu_{R,c}$ is
always non-trivial.

\end{rem}

\begin{rem}

Michael Benedicks pointed to me that a complex atomic measure can indeed
vanish over, say, complex polynomials (for a discussion on this problem see
\cite {BSZ}).

\end{rem}

\appendix

\section{Appendix: infinitesimal perturbations of polynomial maps and hybrid
classes}

The results we obtained in the setting of rational maps have clear
analogous for polynomial maps.  Instead of presenting those straightforward
applications, we will use this Appendix to present a link to
Douady-Hubbard's Theory of polynomial-like maps.

A polynomial-like map of degree $d$
is a holomorphic
degree $d$ ramified covering map $f:U \to U'$, where $U$ and $U'$ are
Jordan disks and $\overline U \subset U'$.  Such maps
were introduced by \cite {D}.
We assume the reader is familiar with this theory.  The filled-in Julia set
of $f$, $K(f)$ is the set of non-escaping points under $f$.

A polynomial of degree $d$ naturally restricts to a polynomial-like map:
just take $U'=\D_r$ for $r$ big enough.  This allows us to consider
polynomial-like perturbations of a polynomial map.
A holomorphic perturbation of a polynomial on a
neighborhood of the filled-in Julia set $K(P)$
of a polynomial $P$ (with, say, $K(P)$ connected)
will usually give rise to a polynomial like map.
Such perturbations are
more restrictive then the ones we consider here.  Nevertheless they
have interesting geometric interpretations.  The reason is that the space
of polynomial-like maps with connected filled-in Julia set fibers over the
space of polynomials with connected Julia set
(modulo conformal equivalence).  Each fiber is called a
hybrid class\footnote {Two polynomial-like maps $f$ and $g$ with connected
Julia set are hybrid equivalent if there is a quasiconformal conjugacy $h$
between $f$ and $g$ on a neighborhood of their filled-in Julia sets
satisfying $\op h|K(f)=0$.}
and the polynomial (or rather its conformal class)
in each hybrid class is called the straightening of the class.

In this case, a polynomial-like perturbation of a polynomial map
may be projected back to the space of
polynomials with the straightening\footnote {Although the
straightening is not in general continuous at polynomial-like maps,
it is continuous at polynomials,
which is enough for our purposes.  This continuity is analogous to the one
of Lemma \ref {qrpert}, and can be proved in a similar way
(see also Remark \ref {qrpertrem}).}.
This approach, while more
geometric, gives slightly weaker
results, since $J(P)=\partial K(P)$, so the space of holomorphic
perturbations on $J(P)$ (which we considered before) is usually bigger.

Let us elaborate a little bit more the case of quadratic
polynomial maps, that is, we will consider polynomials of the form
$p_c=z^2+c$ (see Remark \ref {deg} for higher degree).  To keep things
non-trivial, we will consider only parameters $c$ in the Mandelbrot set
$\MM$, the set of all $c$ such that the critical orbit of $p_c$ is
non-escaping.

According to Lyubich in \cite {universe}, \S 4,
the space $\CC$ of germs\footnote {Two quadratic-like
maps $f$ and $g$ with connected filled-in Julia set give rise to the same
germ if $K(f)=K(g)$ and $f|K(f)=g$.} of normalized
quadratic-like maps with a connected Julia set
has a complex analytic structure.  The tangent space
to a germ $f \in \CC$ is precisely the space of germs of vector
fields which are holomorphic on a
neighborhood of $K(f)$ and whose derivative vanish to order $2$ at $0$.
This space is laminated by
hybrid classes $\HH_c$, $c \in \MM$, which are codimension-one
analytic submanifolds.  The quadratic family
$c \mapsto z^2+c$ is a transverse to this
lamination, see Theorem 4.11 of \cite {universe}.
Finally, the filled-in Julia set of $f \in \CC$ depends
holomorphically along $\HH_c$.  Since the critical orbit of $p_c$ is
contained in $K(p_c)$ we have as before:

\begin{lemma} \label {hyb}

Let $c \in \MM$ and let $v$ be a holomorphic vector field on a
neighborhood of $K(p_c)$.  If $v$ is tangent to the hybrid class of $c$,
say, along an analytic path $f_\lambda \in \HH_c$,
$$
f_0=p_c,\quad\left .\frac {d} {d\lambda} f_\lambda\right |_{\lambda=0}=v,
$$
then $v=\alpha \circ p_c-\alpha Dp_c$ for some qc vector field
$\alpha$, and in
particular, the sequence
$$
\left .\frac {d} {d\lambda} (f_\lambda)^n(0) \right |_{\lambda=0}=
Dp_c^{n-1}(c) \sum_{k=0}^{n-1} \frac {v(p_c^k(c))} {Dp_c^k(p_c(c))}
$$
is uniformly bounded.

\end{lemma}

\begin{rem}

In Section 4.4 of \cite {universe}, it is actually shown
that a tangent vector
field $v$ to the hybrid class of a polynomial-like map $f$ admits a
representation $v=\alpha \circ f-\alpha Df$ satisfying $\op \alpha|K(f)=0$. 
Moreover, those two conditions characterize the tangent space to the hybrid
class of $f$.

\end{rem}

This allows us to characterize the tangent space to the hybrid class of
polynomials satisfying the summability condition:

\begin{lemma} \label {non-tri}

Assume $p_c$ has a summable critical point.  Then $\mu_{p_c,0}$ does not
vanish over the $T_{p_c} \CC$.

\end{lemma}

\begin{pf}

According to Przytycki (see \cite {R}),
under the summability assumption one has that
$K(p_c)$ is a full compact set with empty interior.
By the Mergelyan Approximation Theorem, any
continuous function over $K(p_c)$ can be approximated by a polynomial
in $T_{p_c} \CC$.  The result follows since $\mu_{p_c,0}$ does not vanish
over continuous functions.
\end{pf}

\begin{thm} \label {kernel}

If $p_c$ has a summable critical point then the tangent space to $\HH_c$ is
the kernel of $\mu_{p_c,0}$.

\end{thm}

\begin{pf}

Indeed, $\mu_{p_c,0}$ vanishes over the tangent space to $\HH_c$ by Lemma
\ref {hyb}, but does
not vanish completely by Lemma \ref {non-tri}.
The result follows since $\HH_c$ is codimension-one.
\end{pf}

Let us mention a curious application of transversality of the
quadratic family:

\begin{cor} \label {neq}

Assume the critical point of $p_c$ is summable.  Then
\begin{equation} \label {equa}
\mu_{p_c,c}(1)=\sum_{k=0}^\infty \frac {1} {Dp_c^k(c)} \neq 0.
\end{equation}

\end{cor}

\begin{pf}

Indeed, the tangent to the quadratic family is the set of constant
vector fields.
Transversality implies that the tangent vector field $1$ does not belong to
the tangent space of $\HH_c$.
\end{pf}

(This same result seems to have been independently obtained by Levin using
other techniques.)

In \cite {R}, the same result is obtained with a stronger
condition than summability (and under those conditions, an
interpretation for $\mu_{p_c,c}(1)$
is given in terms of ratios of similarity between Julia and
Mandelbrot sets).

\begin{rem} \label {deg}

The results described here for quadratic-like maps, in particular Theorem
\ref {kernel} and Corollary \ref {neq}
can be also obtained in the context of unicritical
polynomial-like maps of degree $d$ (hybrid conjugate to $z^d+c$).
In this case there is still a lamination of the space $\CC_d$ of germs with
connected Julia set by hybrid classes, and the generalized Mandelbrot set
$\MM_d$ is a global transversal (although a given hybrid class might
intersect $\MM_d$ in more than one point).  In fact for our purposes it is
only necessary to generalize (in a straightforward way)
Sections 4.4 and 4.5 of \cite {universe} which deal
with the description of the tangent space of $\CC_d$.

\end{rem}

\begin{problem}

Let $c \in \MM_d$ be such that $p_c$ is non-hyperbolic.  Is it true that
\begin{equation} \label {equat}
Dp_c^{n-1}(c) \sum_{k=0}^{n-1} \frac {1} {Dp_c^k(c)}
\end{equation}
is unbounded?

\end{problem}

The previous Corollary shows that (\ref {equat}) is unbounded
at least in the summable case.  The work of \cite {ALM} indicates that
this should hold for real quadratic maps whose postcritical set is a
non-minimal set.  By Corollary \ref {struc},
an affirmative answer in the general case
would imply that the Fatou conjecture is valid in this setting.

\end{document}